\documentclass[11pt, twoside]{article}

\usepackage[margin=2.5cm]{geometry}
\usepackage{setspace}

\usepackage{amsmath, amssymb, amsthm}
\usepackage{enumitem}
\usepackage{dsfont}
\usepackage{relsize}

\usepackage{verbatim}
\usepackage{cases}

\usepackage{graphicx,float}
\usepackage[ruled,linesnumbered,nosemicolon,algo2e]{algorithm2e}

\usepackage{hyperref}
\usepackage[svgnames]{xcolor}

\usepackage[numbers,sort,compress]{natbib}

\newtheorem{theorem}{Theorem}

\newtheorem{proposition}[theorem]{Proposition}
\newtheorem{lemma}[theorem]{Lemma}

\theoremstyle{definition}
\newtheorem{assumption}{Assumption}
\newtheorem{definition}{Definition}
\newtheorem{remark}{Remark}

\newcommand{\RR}{\mathbb{R}}

\newcommand{\PP}{\mathbb{P}}
\newcommand{\EE}{\mathbb{E}}

\newcommand{\BB}{\mathcal B}

\newcommand{\norm}[1]{\left\lVert #1 \right\rVert}
\newcommand{\inner}[2]{\left\langle #1, #2 \right\rangle}

\newcommand{\tbo}[2]{\begin{bmatrix}#1 \\ #2\end{bmatrix}}

\newcommand{\argmax}{\mathop{\mathrm{argmax}}}

\newcommand{\tr}{\mathrm{tr}}

\newcommand{\opnorm}[1]{\norm{#1}_{\mathrm{op}}}
\newcommand{\hinfnorm}[1]{\norm{#1}_{\mathcal{H}_\infty}}

\title{Finite-Time Analysis of Supervisory Switching Control via \\Graph-Structured Model Evaluation}
\author{Haoyuan Sun, Ali Jadbabaie\\ \\Massachusetts Institute of Technology}
\date{}

\begin{document}

\maketitle

\begin{abstract}We present a finite-time analysis of supervisory switching control for partially-observed linear dynamical systems.
    The objective is to identify and deploy a stabilizing controller for an unknown plant by switching among $N$ candidate controllers, some of which may destabilize the underlying system.
    While recent non-asymptotic supervisory control methods establish finite-time performance bounds, they only attempt to evaluate the one active controller in each episode, resulting in sample complexity that scales log-linearly with $N$.
    To overcome this bottleneck, we propose a cross-model evaluation approach that simultaneously certifies or falsifies multiple candidate models from a single output trajectory.
    Our criterion has two parts: an observability-based state estimator that isolates the initial condition of each episode from historical states;
    and a short-horizon threshold test on outputs measured along carefully constructed \emph{probing directions} that maintains a sufficient signal-to-noise ratio under unstable dynamics.
    The permissible model-controller pairs that meet the conditions of this criterion form a directed \emph{feedback-information graph}, which we integrate into a graph-structured online learning algorithm to guide switching decisions.
    We prove that our algorithm identifies the matching controller and guarantees finite closed-loop $\mathcal{L}_2$-gain in $O(\alpha \, \mathrm{polylog}(N))$ steps, where $\alpha \le N$ is the graph's independence number, improving upon prior analyses by up to a factor of $N$.
\end{abstract}

\renewcommand{\thefootnote}{\roman{footnote}}
\footnotetext[0]{This work was supported by AFOSR MURI grant \#FA9550-25-1-0375.}
\footnotetext[0]{Corresponding author: Haoyuan Sun (\texttt{haoyuans [at] mit [dot] edu}).}
\renewcommand{\thefootnote}{\arabic{footnote}}

\section{Introduction}
This work presents a non-asymptotic analysis of the supervisory switching control problem.
In many applications involving control of complex systems, such as power systems~\citep{meng2016microgrid}, autonomous vehicles~\citep{aguiar2007trajectory}, and epidemic control~\citep{bin2021hysteresis}, it is impractical for a single controller design to achieve acceptable performance across all operating conditions of the unknown plant.
\textit{Switching control} addresses this challenge by employing a collection of candidate controllers coordinated by a higher-level \textit{supervisor} that monitors system output measurements and dynamically switches to the controller best suited to the current system behavior~\citep{liberzon2003switching}.
A well-studied example of this principle is the \textit{estimator-based supervision}~\citep{morse1996supervisory,morse1997supervisory}, where candidate controllers are designed for a family of models that cover the plant's uncertainty set.
In practice, this covering is typically finite to balance between approximation accuracy and computational cost~\citep{anderson2000multiple}.
A \textit{multi-estimator} then simulates the predicted responses of all candidate models in parallel and compares them against the plant outputs. 
Subject to a \textit{dwell time}, the supervisor from time to time switches to the controller whose model is most consistent with the observations.
However, the classical analysis only establishes asymptotic stability guarantees, and it does not precisely characterize when these estimator-based evaluations succeed, nor how long it takes to find a suitable controller.
 
In this work, we introduce a cross-model criterion that evaluates multiple candidates from a single uninterrupted trajectory.
We derive finite-time guarantees for this criterion even when the active controller destabilizes the underlying plant.
We organize the model-controller pairs that admit valid simultaneous evaluations into a directed \emph{feedback-information graph}.
By exploiting this graph structure, our online algorithm provably identifies the matching stabilizing controller without having to deploy each candidate controller to the feedback-loop.
 
Our analysis is inspired by the extensive literature in learning-based control of unknown dynamical systems.
Within this space, one approach builds on a non-asymptotic theory of system identification from a single trajectory~\citep{simchowitz2018learning,oymak2019non} and derives control performance bounds via the certainty equivalence principle~\citep{mania2019certainty,simchowitz2020naive}.
A second perspective frames the control problem in terms of online learning and directly optimizes the decision sequences~\citep{abbasi2011regret,agarwal2019online,muehlebach2026the}.
However, these analyses break down in supervisory switching control because they rely on closed-loop stability or independence across control actions.
Such assumptions are typically enforced by assuming a known stabilizing controller or, in the case of episodic reinforcement learning (RL)~\citep{jin2018q,du2021bilinear}, by periodically resetting the environment.
But in supervisory switching control, candidate controllers can destabilize the unknown plant, and the learner must be able to handle potentially destabilizing switching decisions along a single trajectory without resets.

Several recent works have attempted to tackle this challenge under instability, for example by injecting exponentially large probing signals~\citep{chen2021black} or by designing stability certificates on system states~\citep{li2023online,kim2025online}.
However, these techniques require full state observation and cannot accommodate our setting with partially-observed states.
Our prior work~\citep{sun2024least,sun2026online} proposed solutions in the partially-observed setting by employing two evaluation criteria, where the first detects destabilizing controllers, and the second applies least-squares system identification to the remaining candidates.
By leveraging system observability to cancel out the influence of earlier switching decisions, our prior work derived a bandit-style selection rule that identifies the matching controller in $O(N \, \mathrm{polylog} N)$ steps~\citep{sun2026online}.
Crucially, however, both of our previous algorithms evaluate only a single candidate model at a time and discard whatever information the trajectory carries regarding the remaining models.
To address this shortcoming, we now propose a novel evaluation criterion that tests the observed trajectory against multiple candidate models at once, directly extending the cross-model evaluation in a classical multi-estimator to a provable finite-time framework.
Inspired by bandits with side observations~\citep{mannor2011bandits}, our online algorithm exploits the outcomes of these simultaneous evaluations to achieve improved sample complexity bounds.
Additionally, we defer a more extensive survey of related literature to Appendix~\ref{sec:related-work}.

\paragraph{Contributions.}
We develop an online learning framework for supervisory switching control that establishes non-asymptotic sample complexity and stability guarantees for partially observed linear systems with a collection of $N$ candidate controllers.
Operating along a single trajectory without interruptions or resets, our approach simultaneously certifies or falsifies multiple candidate models for each episode with predetermined length --- even when the active controller destabilizes the underlying plant.
By capturing these permissible cross-model evaluations as a directed \emph{feedback-information graph} over the $N$ candidates, our online switching algorithm identifies the matching controller and guarantees finite $\mathcal{L}_2$-gain stability within $O(\alpha \, \mathrm{polylog}(N))$ time steps, where $\alpha \le N$ is the graph's independence number.
Specifically, our main technical contributions are three-fold:
\begin{enumerate}[wide, labelindent=3pt, label=\arabic*), itemsep=2pt]
    \item In Section~\ref{sec:criteria}, we design a two-stage evaluation criterion (Algorithm~\ref{alg:moment-threshold}) that validates multiple candidate models $i$ from a single output trajectory generated under an active controller $j$ (where $i$ need not equal $j$), generalizing our prior work~\citep{sun2026online} which only evaluates the case where $i = j$.
    First, we leverage system observability and a state observer to reconstruct the system state and subtract out influences from prior switching decisions.
    Second, we project the outputs onto carefully constructed \emph{probing directions} that preserve a high signal-to-noise ratio in response to an exploratory signal over a short horizon.
    Such probing directions enable a moment-based threshold test that certifies or rejects candidate models with high probability (Proposition~\ref{thm:moment-threshold}).

    \item In Section~\ref{sec:main-result}, we characterize the permissible model-controller pairs satisfying our observability and probing conditions via a directed feedback-information graph.
    Combining this graph structure with an importance-weighted multiplicative weights update (Algorithm~\ref{alg:meta-graph}), we dynamically balance exploration and exploitation across candidate controllers.
    We prove that our algorithm identifies the matching controller in $\mathcal{O}(\alpha \log N)$ episodes while maintaining finite closed-loop $\mathcal{L}_2$-gain (Theorem~\ref{thm:graph-exp3}).
    This improves upon the sample complexity of our previous work~\citep{sun2026online} by a factor of $N$ when the feedback-information graph is complete and gracefully recovers our prior $O(N \log N)$ guarantee in the worst case where each controller can only evaluate its own model.

    \item Finally, in Appendix~\ref{app:simulation}, we illustrate the efficacy of Algorithm~\ref{alg:meta-graph} on a benchmark from~\citet{anderson2000multiple}.
    There, the feedback-information graph is complete (i.e., $\alpha = 1$), and Algorithm~\ref{alg:meta-graph} yields a significantly smaller transient than both classical estimator-based supervision and our prior algorithm that only evaluates one candidate at a time~\citep{sun2026online}.
\end{enumerate}

\section{Problem settings}

In this work, we consider a collection of discrete-time, partially-observed linear systems $\left\{(C_i, A_i, B_i)\right\}_{i=1}^N$.
The unknown true system parameters $(C_\star, A_\star, B_\star)$ belong to this collection, and have index $i_\star$. 
Each of these systems evolves according to the following discrete-time dynamics:
\begin{align*}
    x_{t+1} = A_i x_t + B_i \bar{u}_t + w_t, \qquad y_t = C_i x_t + \eta_t,
\end{align*}
where the dimensions are $x_t \in \RR^{d_x}, \bar{u}_t \in \RR^{d_u}$ and $y_t \in \RR^{d_y}$.
We assume zero initialization $x_0 = 0$ and independent Gaussian noise for both the process disturbance $w_t \sim \mathcal{N}(0, \sigma_w^2 I_{d_x \times d_x})$, and observation noise $\eta_t \sim \mathcal{N}(0, \sigma_\eta^2 I_{d_y \times d_y})$.

Each candidate model is paired with a stabilizing candidate linear controller that provides satisfactory performance for the associated model.
However, a mismatched pair of plant model and controller may result in an unstable closed-loop system.

The supervisor switches between the candidate controllers, and we model this behavior with the \textit{multi-controller} architecture $K(p_t; \check{x}_t, y_t)$ in \citep{hespanha2001tutorial}, where $\check{x}_t$ is the internal state of the controller, $p_t \in [N]$ is a piece-wise constant switching signal indicating the active controller, and $y_t$ is the system's output.
And to ensure persistency of excitation, we inject an independent, additive Gaussian exploratory signal $u_t \sim \mathcal{N}(0, \sigma_u^2 I_{d_u \times d_u})$, so that the total control signal into the open-loop system is $\bar{u}_t = u_t + K(p_t, \check{x}_t, y_t)$.

\begin{figure}[!htb]
    \centering
    \includegraphics[scale=0.9]{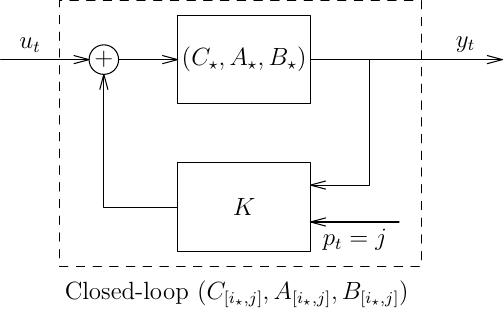}
    \vspace{-1em}
    \caption{Illustration of a linear switching system (with a fixed switching signal)}
    \label{fig:switching-loop}
\end{figure}

For every pair $1 \le i, j \le N$, we denote $M_{[i, j]} = \left(C_{[i, j]}, A_{[i, j]}, B_{[i, j]}\right)$ as the closed-loop system resulting from applying the $j$th candidate controller to the $i$th open-loop model.
As illustrated by Figure~\ref{fig:switching-loop}, applying a constant switching signal $p_t = j$ on the switching system results in a closed-loop modeled by $M_{[i_\star, j]}$.
Furthermore, we can pre-compute every possible closed-loop system $M_{[i,j]}$ from the collection of candidate models and controllers.
Consequently, under a piecewise-constant switching signal, our observations are generated from one of these $N^2$ possible processes.

To quantitatively analyze the finite-time performance of this switching system, we need \textit{observability} of the closed-loop systems.

\begin{definition}
\label{assum:observability}
    We say that a linear system $(C, A, B)$ is $\varepsilon_c$-strictly observable with index $\nu$ if the observability matrix satisfies
    $ \sigma_{\min}\left(\left[CA^{\nu-1}; \dots; CA; C\right]\right) \ge \varepsilon_c.$
\end{definition}

Following the conventions in the learning-for-control literature (e.g. \cite{chen2021black, sarkar2021finite}), our assumptions are non-asymptotic counterparts to the standard control definitions.
Under the premise that each candidate model's associated controller results in satisfactory performance, we assume that each closed-loop system with matched model-controller $M_{[j,j]}$ is observable and stable.
Furthermore, we assume that we can distinguish the matched closed-loop $M_{[j,j]}$'s with a sufficiently large margin.

\begin{assumption}
\label{assum:matched-controller}
    For each $j \in [N]$, the closed-loop system $M_{[j, j]}$ is $\varepsilon_c$-strictly observable and stable so that its system $\mathcal{H}_\infty$-norm is bounded by $\Lambda$.
    Also, there exists $h, \gamma > 0$ so that for each $j$ and $i \neq j$, we have $\opnorm{G_{[i,j]} - G_{[j,j]}} > \gamma$, where $G := \left[C A^{h-1} B, \dots, C A B, C B\right]$ is the response matrix to inputs of window $h$.
\end{assumption}
\vspace{-0.5em}
Additional notations and technical lemmas can be found in Appendix~\ref{sec:prelim}.
 
\section{Cross-model evaluation criterion}
\label{sec:criteria}
In this section, we present the core technical innovation of our online algorithm for supervisory switching control: a cross-model evaluation criterion (Algorithm~\ref{alg:moment-threshold}) that provably certifies or falsifies multiple candidate models from a single trajectory.
Our online algorithm (Algorithm~\ref{alg:meta-graph}) employs this criterion in episodes of fixed length $\tau$, a constant specified in Proposition~\ref{thm:moment-threshold} to ensure that sufficient measurements are collected before each switching decision.
During each episode, an active candidate controller is applied to the plant, and the resulting output trajectory is used to evaluate all candidate models that admit a valid \emph{probing direction}.
Unlike our previous work~\citep{sun2026online}, the evaluated models need not correspond to the active controller, thereby enabling simultaneous multi-model feedback per episode.
These permissible model-controller pairs form a directed \emph{feedback-information graph}, which guides switching decisions at the end of each episode.

We recall that existing learning-based control techniques do not directly apply to our setting because they typically assume system stability or independence of control actions over time.
Our evaluation criterion overcomes these challenges through a two-pronged approach.
First, we estimate each episode's initial state from the output trajectory, allowing us to isolate the effects of past episodes.
Second, we apply a moment-based threshold test along carefully selected probing directions over a short time horizon, ensuring that model evaluation succeeds with high probability even when the closed-loop system is unstable.
Specifically, when the evaluated model $i$ matches the true plant model $i_\star$, the state-estimation step cancels out the influence of historical states, and the moment-based test verifies that the candidate model is consistent with the observations.
Conversely, when the evaluated model $i$ does not match $i_\star$, the probing directions align with the discrepancy in model dynamics over the short horizon, so that the moment-based thresholds are exceeded regardless of the historical state trajectory.

\subsection{Part 1: State estimation via observability}
\label{sec:state-estimation}
In the first stage of our criterion, we seek to cancel out the effects of past episodes by reconstructing the system state from output measurements.
For convenience, denote the current episode's outputs as $y_1, \dots, y_\tau$, and let $i$ and $j$ index the evaluated candidate model and the active controller, respectively.
Because accurately estimating the state is generally intractable without knowing the true plant dynamics, we design our estimator using the dynamics of the candidate pair $[i, j]$. 
This way, we only guarantee a small state-estimation error when the evaluated model matches the ground truth ($i = i_\star$), leaving the detection of mismatched models ($i \neq i_\star$) to the second stage of our criterion.

Under the true closed-loop system $M_{[i_\star, j]}$, expanding the dynamics over the first $\nu$ steps yields
\begin{equation}
\label{equ:obs-dynamics}
\begin{bmatrix} y_\nu \\ \vdots \\ y_2 \\ y_1 \end{bmatrix}
=
\mathcal{O}_{[i_\star, j]} x_1
+
\tbo{\mathcal{T}^{(\nu-1)}_{[i_\star, j]}}{\mathbf{0}_{1 \times (\nu-1)}}
\begin{bmatrix} w_{\nu-1} \\ \vdots \\ w_2 \\ w_1  \end{bmatrix}
+
\tbo{\mathcal{T}^{(\nu-1)}_{[i_\star, j]}}{\mathbf{0}_{1 \times (\nu-1)}}
\begin{bmatrix} B_{[i_\star, j]}u_{\nu-1} \\ \vdots \\ B_{[i_\star, j]}u_2 \\ B_{[i_\star, j]}u_1  \end{bmatrix}
+
\begin{bmatrix} \eta_\nu \\ \vdots \\ \eta_2 \\ \eta_1  \end{bmatrix},
\end{equation}
where $\mathcal{O}_{[i, j]} := [C_{[i, j]}A_{[i, j]}^{\nu-1}; \dots; C_{[i, j]}A_{[i, j]}; C_{[i, j]}]$ denotes the $\nu$-step observability matrix and $\mathcal{T}^{(k)}_{[i, j]}$ is the corresponding $k \times k$ Toeplitz matrix.
For simplicity, let $\mathcal{T}^{(\nu-1)}_{[i, j]}[u]$ denote the third term on the right-hand side of~\eqref{equ:obs-dynamics}, which captures the response to the exploratory inputs $u_1, \dots, u_{\nu-1}$.
When $i = i_\star$, multiplying both sides of~\eqref{equ:obs-dynamics} by the Moore-Penrose pseudo-inverse $\mathcal{O}_{[i, j]}^\dagger$ yields an estimate of the initial state:
\begin{equation}
\label{equ:obs-dynamics-short}
    \hat{x}_1 := \mathcal{O}_{[i, j]}^\dagger \left([y_\nu; \dots; y_1] - \mathcal{T}^{(\nu-1)}_{[i, j]}[u]\right) = x_1 + \mathcal{O}_{[i, j]}^\dagger \xi,
\end{equation}
where $\xi$ aggregates the process noise $w$ and measurement noise $\eta$ over the $\nu$-step window.

When the condition in Definition~\ref{assum:observability} holds, the closed-loop system $M_{[i, j]}$ is $\varepsilon_c$-strictly observable, which means that $\opnorm{\mathcal{O}_{[i, j]}^\dagger} \le 1/\varepsilon_c$.
However, if $M_{[i, j]}$ is unstable, the operator norm of the Toeplitz matrix $\mathcal{T}^{(\nu-1)}_{[i, j]}$ grows exponentially with $\nu$, which can amplify the noise term $\mathcal{O}_{[i, j]}^\dagger \xi$ in~\eqref{equ:obs-dynamics-short}.
To dampen this amplified noise, we can design a \emph{state observer} gain $L_{[i, j]}$ such that $A_{[i, j]} - L_{[i, j]}C_{[i, j]}$ is exponentially stable, and propagate the state estimate forward:
\begin{align*}
    \hat{x}_{t+1} &= (A_{[i,j]} - L_{[i,j]} C_{[i,j]}) \hat{x}_t + B_{[i,j]} u_t + L_{[i,j]} y_t.
\end{align*}
When $i = i_\star$, the estimation error $x_t - \hat{x}_t$ evolves according to the stable dynamics
\begin{align*}
    x_{t+1} - \hat{x}_{t+1} &= (A_{[i,j]} - L_{[i,j]} C_{[i,j]}) (x_t - \hat{x}_t) + w_t - L_{[i,j]} \eta_t,
\end{align*}
which expands after $t$ steps to
\[
x_t - \hat{x}_t = (A_{[i,j]} - L_{[i,j]} C_{[i,j]})^{t-1} \mathcal{O}_{[i,j]}^\dagger \xi + \sum_{s=1}^{t-1} (A_{[i,j]} - L_{[i,j]} C_{[i,j]})^{t-s-1} (w_s - L_{[i,j]} \eta_s).
\]
Because the initial estimation error $\mathcal{O}_{[i,j]}^\dagger \xi$ decays exponentially fast under $(A_{[i,j]} - L_{[i,j]} C_{[i,j]})^{t-1}$, the state-estimation error becomes $\sigma_r^2$-sub-Gaussian after $\tau_1 = \mathcal{O}(\nu)$ burn-in steps, where $\sigma_r^2$ depends on the norms of $A_{[i,j]}$ and $A_{[i,j]} - L_{[i,j]} C_{[i,j]}$.
Hence, we can produce an accurate estimate on the system state when $i = i_\star$.

\hspace{-1em}
\subsection{Part 2: Moment-based threshold}
\label{sec:moments}
In the second part of our criterion, we determine whether the evaluated model $i$ matches the true plant $i_\star$ by inspecting how the exploratory inputs $u$ influence the output trajectory over a short horizon of $h$ steps.
Specifically, we compare the actual observed outputs against the outputs predicted by candidate model $i$ initialized at the estimated state $\hat{x}_{t-h}$.
When $i = i_\star$, the residual differences are purely noise and independent of the exploratory inputs $u$.
But when $i \neq i_\star$, the residuals would retain a correlation with $u$.

To formalize this idea, we unroll the true closed-loop dynamics $M_{[i_\star, j]}$ over a horizon of $h$ steps:
\hspace{-0.5em}
\begin{align*}
    y_t &= C_{[i_\star, j]} A_{[i_\star, j]}^{h} x_{t-h} \! + \! \underbrace{\sum_{s=1}^h C_{[i_\star, j]} A_{[i_\star, j]}^{s-1} B_{[i_\star, j]} u_{t-s}}_{G_{[i_\star, j]} z_t} \! + \! \underbrace{\sum_{s=1}^h C_{[i_\star, j]} A_{[i_\star, j]}^{s-1} w_{t-s} + \! \eta_t}_{e_t} \\
    &:= C_{[i_\star, j]} A_{[i_\star, j]}^{h} x_{t-h} + G_{[i_\star, j]} z_t + e_t,
\end{align*}
where $z_t := [u_{t-h}; \dots; u_{t-1}]$ is a sliding window of exploratory inputs, $e_t$ aggregates the process and measurement noise, and $G_{[i_\star, j]}$ denotes the $h$-step Markov parameter matrix:
\[G_{[i_\star, j]} := \left[C_{[i_\star, j]} A^{h-1}_{[i_\star, j]} B_{[i_\star, j]}, \dots, C_{[i_\star, j]} A_{[i_\star, j]} B_{[i_\star, j]}, C_{[i_\star, j]} B_{[i_\star, j]}\right].\]
Given a state estimate $\hat{x}_{t-h}$, we construct the \emph{residual output trajectory} $\widetilde{y}_t$ by subtracting the candidate model's predicted rollout plus nominal input response:
\begin{align*}
    \widetilde{y}_t 
    &:= y_t - C_{[i, j]} A_{[i, j]}^{h}\hat{x}_{t-h} - G_{[i, j]} z_t \\
    &= (G_{[i_\star, j]} - G_{[i, j]}) z_t + \left(C_{[i_\star, j]} A_{[i_\star, j]}^{h} x_{t-h} - C_{[i, j]} A_{[i, j]}^{h}\hat{x}_{t-h}\right) + e_t \\
    &:= \Delta_{i_\star, i} z_t + \zeta_t + e_t,
\end{align*}
where $\Delta_{i_\star, i} := G_{[i_\star, j]} - G_{[i, j]}$ is the Markov parameter gap and $\zeta_t$ captures the state-prediction discrepancy at the start of the $h$-step window.
Crucially, because $z_t$ consists of fresh exploratory inputs generated from time $t-h$ to $t-1$, both $\zeta_t$ and $e_t$ are independent of $z_t$.
We then inspect $\widetilde{y}_t$ under two scenarios where the evaluated model matches the true system ($i = i_\star$), or does not match the true system ($i \neq i_\star$).

\noindent \textbf{Case 1: $i = i_\star$.}
Here $\Delta_{i_\star, i} = 0$ and $\zeta_t = C_{[i, j]} A_{[i, j]}^{h}(x_{t-h} - \hat{x}_{t-h})$.
Because $A_{[i,j]}$ may be unstable, $C_{[i, j]} A_{[i, j]}^{h}$ can amplify the estimation error $x_{t-h} - \hat{x}_{t-h}$.
Thus, we project $\widetilde{y}_t$ along a unit probing direction $a$ that avoids the explosive modes of $A_{[i,j]}$, satisfying $\sum_{s=0}^{h} \norm{a^\top C_{[i, j]} A_{[i, j]}^{s}}^2 \le \Lambda^2$.
Then, for $t \ge \tau_1 + h$, $a^\top \widetilde{y}_t$ is zero-mean and sub-Gaussian with variance at most $(\sigma_r^2 + \sigma_w^2)\Lambda^2 + \sigma_\eta^2$.
Consequently, for any unit vector $b$, the cross-term $(a^\top \widetilde{y}_t)(b^\top z_t)$ is zero-mean and sub-exponential, allowing us to upper-bound both the empirical second moment of $a^\top \widetilde{y}_t$ and the empirical cross-moment of $(a^\top \widetilde{y}_t)(b^\top z_t)$ with high probability.

\noindent \textbf{Case 2: $i \neq i_\star$.}
In this case, the Markov parameter gap $\Delta_{i_\star, i}$ is nonzero.
Then there exists a unit probing direction $a_{i_\star}$ along which $\norm{a_{i_\star}^\top \Delta_{i_\star, i}} > \gamma$ for some positive margin $\gamma$.
We set probing direction $b_{i_\star} = \Delta_{i_\star, i}^\top a_{i_\star} / \norm{a_{i_\star}^\top \Delta_{i_\star, i}}$ on the exploratory signal to align with $a_{i_\star}$.
Using the independence of $z_t \sim \mathcal{N}(0, \sigma_u^2 I)$ from $\zeta_t$ and $e_t$, the population cross-moment along $(a_{i_\star}, b_{i_\star})$ satisfies
\[\EE[(a_{i_\star}^\top \widetilde{y}_t)(b_{i_\star}^\top z_t)] = \EE[a_{i_\star}^\top \Delta_{i_\star, i} z_t z_t^\top b_{i_\star}] = \sigma_u^2 \norm{a_{i_\star}^\top \Delta_{i_\star, i}} > \sigma_u^2 \gamma.\]
If $\zeta_t$ diverges under an unstable closed-loop dynamics, the empirical second moment of $a_{i_\star}^\top \widetilde{y}_t$ will exceed our variance threshold.
Otherwise, $a_{i_\star}^\top \zeta_t$ remains bounded, ensuring the empirical cross-moment concentrates above $\sigma_u^2 \gamma / 2$.
These requirements can be formalized in Definition~\ref{def:probing}.

\begin{definition}[Probing Directions]
\label{def:probing}
    A model-controller pair $(i, j)$ admits \emph{($\gamma, \Lambda$)-probing directions} if, for every alternative model $i' \neq i$, there exist unit vectors $a_{i'}$ and $b_{i'}$ satisfying:
    \begin{itemize}[wide, labelindent=3pt, itemsep=2pt]
    \item $\norm{a_{i'}^\top \Delta_{i', i}} > \gamma$ and $b_{i'} = \Delta_{i', i}^\top a_{i'} / \norm{a_{i'}^\top \Delta_{i', i}}$, where $\Delta_{i', i} := G_{[i', j]} - G_{[i, j]}$.
    \item $\sum_{s=0}^{h} \norm{a_{i'}^\top C_{[i, j]} A_{[i, j]}^{s}}^2 \le \Lambda^2 \quad \forall i' \neq i$.
    \end{itemize}
\end{definition}

\begin{remark}
Definition~\ref{def:probing} formalizes the sufficient conditions for ensuring strong signal-to-noise ratio in cross-model evaluations:
\begin{itemize}[wide, labelindent=3pt, itemsep=2pt]
    \item The vectors $a_{i'}, b_{i'}$ respectively project the outputs and the exploratory signal to expose the difference between the $h$-step Markov parameters of candidate models $i$ and $i'$, ensuring a separation margin of at least $\gamma$ along this direction.
    \item When $i \neq j$, the closed-loop matrix $A_{[i, j]}$ may be unstable. Steering $a_{i'}$ away from the unstable modes of $A_{[i, j]}$ prevents an $h$-step rollout from amplifying the state-estimation error $x_{t-h} - \hat{x}_{t-h}$ and process noise $e_t$ beyond a factor of $\Lambda$, ensuring that the discrepancy signal is not overwhelmed.
\end{itemize}
The second condition is mild because it only restricts growth over a short horizon $h \ll \tau$, rather than directly demanding stability along the direction $a_{i'}$ for all $\tau$ steps.
Moreover, for self-evaluation ($i = j$), the nominal closed-loop system $M_{[j, j]}$ is stable under Assumption~\ref{assum:matched-controller} and automatically meets the second condition for any unit vector $a_{i'}$.
As a result, every controller $j$ can always evaluate its own candidate model $i = j$.
\end{remark}

When evaluated model $i$ and active controller $j$ admit valid probing directions, the computed residuals are expected to have the following behaviors:
\begin{itemize}[wide, labelindent=3pt, itemsep=2pt]
    \item If $i = i_\star$, then for all $i' \neq i$, the population moments satisfy $\EE[(a_{i'}^\top \widetilde{y}_t)^2] \le (\sigma_r^2 + \sigma_w^2)\Lambda^2 + \sigma_\eta^2$ and $\EE[(a_{i'}^\top \widetilde{y}_t)(b_{i'}^\top z_t)] = 0$.
    \item If $i \neq i_\star$, then for the true index $i' = i_\star$, either the second moment $\EE[(a_{i_\star}^\top \widetilde{y}_t)^2]$ is large, or $\EE[(a_{i_\star}^\top \widetilde{y}_t)(b_{i_\star}^\top z_t)] > \sigma_u^2 \gamma$.
\end{itemize}
The state-estimation and moment-based thresholds together make up our evaluation criterion, which is stated formally in the next section.

\subsection{Criterion's theoretical guarantees}
Combining the state estimation and empirical moments yields our cross-model evaluation criterion (Algorithm~\ref{alg:moment-threshold}).
We can state its guarantees subject to observability (Definition~\ref{assum:observability}) and the existence of probing directions (Definition~\ref{def:probing}).
The complete proof of this result is deferred to Appendix~\ref{sec:proof-criteria}.

\begin{algorithm2e}
\KwIn{Evaluated model $i$, active controller $j$, confidence $\delta_1$, trajectory $y_{1:\tau}, u_{1:\tau}$ with $\tau = h+\tau_1+\tau_2$.}
\For{$t = h+\tau_1+1, \dots, h+\tau_1+\tau_2$}{
    Let $s_0 = t - h - \tau_1$ and initialize
    \begin{equation}
    \label{equ:obs-dynamics-alg}
    \hat{x}_{s_0} := \mathcal{O}_{[i, j]}^\dagger \left([y_{s_0+\nu-1}; \dots; y_{s_0}] - \left[\mathcal{T}^{(\nu-1)}_{[i, j]}; \mathbf{0}\right] [B_{[i, j]}u_{s_0+\nu-2}; \dots; B_{[i, j]}u_{s_0}]\right)
    \end{equation}
    and roll forward
    $
        \hat{x}_{s+1} = (A_{[i,j]} - L_{[i,j]} C_{[i,j]}) \hat{x}_s + B_{[i,j]} u_s + L_{[i,j]} y_s
    $
    to obtain $\hat{x}_{t-h}$.\;
    Compute $\widetilde{y}_t = y_t - C_{[i, j]} A_{[i, j]}^{h}\hat{x}_{t-h} - G_{[i, j]} z_t$.\;
    Measure along the probing directions with $r_{i', t} = a_{i'}^\top \widetilde{y}_t$, $\bar{z}_{i', t} = b_{i'}^\top z_t$ for all $i' \neq i$.\;
}
\ForAll{$i' \neq i$}{
    Compute $\Phi^{(1)}_{i'} = \frac{1}{\tau_2} \sum_{t=h+\tau_1+1}^{\tau} r_{i', t}^2$ and $\Phi^{(2)}_{i'} = \frac{1}{\tau_2} \sum_{t=h+\tau_1+1}^{\tau} r_{i', t} \bar{z}_{i', t}$\;
    \If{$\Phi^{(1)}_{i'} > 4(\Lambda^2 (\sigma_r^2 +\sigma_w^2) + \sigma_\eta^2) \log(2\sqrt{2} N /\delta_1)$ \textbf{or} $\Phi^{(2)}_{i'} > \sigma_u^2 \gamma / 2$}{\Return $0$ \tcp*{Reject model $i$}}
}
\Return $1$ \tcp*{Certify model $i$}
\caption{Cross-Model Evaluation Criterion $S(y_{1:\tau}; i, j)$}
\label{alg:moment-threshold}
\end{algorithm2e}

\begin{proposition}
\label{thm:moment-threshold}
    Suppose the conditions in Definitions~\ref{assum:observability} and~\ref{def:probing} hold for evaluated model $i$ and active controller $j$. Conditioned on any state history and initial condition, for an episode of length $\tau = h + \tau_1 + \tau_2$ with sample horizon
    \[\tau_2 \asymp \frac{\Lambda^2(\sigma_r^2 + \sigma_w^2) + \sigma_\eta^2}{\sigma_u^2 \gamma^2} h \log(N/\delta_1)^2,\]
    the evaluation criterion $S(y_{1:\tau}; i, j)$ in Algorithm~\ref{alg:moment-threshold} satisfies:
    \begin{itemize}[wide, labelindent=3pt]
    \item In the matched case where $i = i_\star$, then $S(y_{1:\tau}; i, j) = 1$ with probability at least $1 - \delta_1$.
    \item In the mismatched case where $i \neq i_\star$, then $S(y_{1:\tau}; i, j) = 0$ with probability at least $1 - \delta_1$.
    \end{itemize}
\end{proposition}

In Proposition~\ref{thm:moment-threshold}, the factors $\frac{\sigma_r^2 + \sigma_w^2}{\sigma_u^2}$ and $\frac{\sigma_\eta^2}{\sigma_u^2}$ capture the effective inverse signal-to-noise ratio of the exploratory input $u_t$ relative to the state-estimation error and system noise, while the $1/\gamma^2$ scaling matches the standard sample complexity for parameter estimation in statistical learning.

\section{Main algorithm and performance guarantees}
\label{sec:main-result}
In the previous section, we established the cross-model evaluation criterion (Algorithm~\ref{alg:moment-threshold}) and bounded its performance with Proposition~\ref{thm:moment-threshold}.
In particular, when a model-controller pair $(i, j)$ satisfies the conditions in Definitions~\ref{assum:observability} and~\ref{def:probing}, deploying controller $j$ allows us to reliably test whether candidate model $i$ matches the ground truth $i_\star$.
These permissible pairs naturally induce a directed graph over the candidate collection, which we formalize below.

\begin{definition}[Feedback-Information Graph]
\label{def:graph}
    Given a collection of $N$ candidate models and controllers indexed by $\{1, \dots, N\}$, the \emph{feedback-information graph} is a directed graph $G = ([N], E)$ in which a directed edge $(j \to i) \in E$ exists if and only if the model-controller pair $(i, j)$ satisfies the conditions in Definitions~\ref{assum:observability} and~\ref{def:probing}.
\end{definition}

Intuitively, an edge $j \to i$ indicates that applying controller $j$ generates trajectory data capable of validating or falsifying model $i$.
Under Assumption~\ref{assum:matched-controller}, $M_{[j, j]}$'s are both observable and stable, which means that the feedback-information graph $G$ contains all self-loops $j \to j$ for $j \in [N]$.

Given the feedback-information graph $G = ([N], E)$, we now present our online supervisory switching algorithm (Algorithm~\ref{alg:meta-graph}).
Inspired by the EXP3-IX algorithm for bandits with side observations~\citep{kocak2014efficient,lattimore2019bandit}, our algorithm maintains a weight $w_{\ell, i}$ for each candidate model $i \in [N]$ and applies a multiplicative weights update procedure at the end of each episode $\ell$ using the criterion's outcomes.
To account for the asymmetric feedback structure of $G$, we normalize each model's loss estimate by its total observation probability $q_{\ell, i} = \sum_{j:\,(j \to i) \in \mathcal{E}} p_{\ell, j}$, regularized by an implicit exploration parameter $\beta > 0$.

We emphasize that standard EXP3 analysis do not directly apply in our setting due to the coupling of switching decisions across episodes and the potential instability of closed-loop systems.
To address this challenge, our cross-model evaluation criterion is designed to have high probability of success regardless of the episode's initial state (Proposition~\ref{thm:moment-threshold}).
So, our algorithm can isolate per-episode loss estimates that are not influenced by past switching decisions.

\begin{algorithm2e}
\DontPrintSemicolon
\KwIn{Candidate set $[N]$, feedback-information graph $G = ([N], E)$, parameters $\beta, \eta > 0$.}
Initialize weights $w_{1, i} = 1$ for all $i \in [N]$\;
\For{episode $\ell = 1, 2, \dots, L$}{
    Compute the controller sampling probabilities and model observation probabilities:
    \[
      p_{\ell,i} := \frac{w_{\ell,i}}{\sum_{k=1}^N w_{\ell,k}},
      \qquad
      q_{\ell,i} := \sum_{j:\,(j\to i) \in \mathcal{E}} p_{\ell,j} \quad \forall i \in [N]
    \]\vspace{-1.5\baselineskip}\;
    Sample active controller $j_\ell \sim p_\ell$ and apply it for $\tau$ steps to collect output trajectory $\vec{y}^{(\ell)} := y_{\tau(\ell-1)+1 : \tau \ell}$\;
    \ForAll{$i \in [N]$}{
        Compute the importance-weighted loss estimate and update weights:
        \[
          \hat{Z}_{\ell, i} := \frac{1 - S(\vec{y}^{(\ell)}; i, j_\ell)}{q_{\ell, i} + \beta} \mathbb{I}\big\{(j_\ell \to i) \in E\big\},
          \qquad
          w_{\ell+1, i} := w_{\ell, i} \exp(-\eta \hat{Z}_{\ell, i})
        \]\vspace{-1.5\baselineskip}\;
    }
}
Commit to the most frequently selected controller $\hat{i} \in \argmax_{i \in [N]} \sum_{\ell=1}^L \mathbb{I}\{j_\ell = i\}$\;
\caption{Graph-structured algorithm for supervisory switching control}
\label{alg:meta-graph}
\end{algorithm2e}

The sample complexity of Algorithm~\ref{alg:meta-graph} scales with the \emph{independence number} $\alpha(G)$ of the feedback-information graph, defined as the cardinality of the largest vertex subset $I \subseteq [N]$ such that no two distinct vertices in $I$ are connected by a directed edge in $E$.

\begin{theorem}
\label{thm:graph-exp3}
    Let $\alpha := \alpha(G)$ be the independence number of $G$, and choose $\bar{\alpha} \ge \alpha, \beta = 1/(80\bar\alpha)$ and $\eta = \beta/2$.
    Suppose each episode has length $\tau$ given from Proposition~\ref{thm:moment-threshold}, and the total number of episodes is $L = O(\bar\alpha \log (N/\delta_2))$.
    Then, with probability at least $1-\delta_2$, Algorithm~\ref{alg:meta-graph} satisfies:
    \begin{itemize}[wide, labelindent=3pt]
    \item After $L$ episodes, the algorithm commits to the matching controller, $\hat{i} = i_\star$.
    \item There exist constants $C_0, C_1 > 0$ such that for any horizon $T > \tau L$,
    \begin{align*} 
        \sum_{t=1}^T \norm{x_t}^2 
        \le C_0 \sum_{t=1}^{\tau L} \big(\norm{u_t}^2 + \norm{w_t}^2\big)
        + C_1 \sum_{t=\tau L + 1}^T \big(\norm{u_t}^2 + \norm{w_t}^2\big).
    \end{align*}
    \end{itemize}
\end{theorem}

\begin{proof}
    We outline the key steps in the proof and defer the complete version to Appendix~\ref{sec:proof-graph-bandit}.
    Since computing the independence number is NP-hard in general, our analysis allows any upper bound $\bar{\alpha} \ge \alpha$, which is easier to obtain.
    
    We first apply a standard multiplicative-weights telescoping argument by tracking the potential function $W_\ell := \sum_{k=1}^N w_{\ell, k}$.
    This argument yields an upper bound on $\sum_{\ell=1}^L A_\ell - \sum_{\ell=1}^L \hat{Z}_{\ell, i_\star}$, where $A_\ell := \sum_{k=1}^N p_{\ell,k} \hat{Z}_{\ell, k}$ denotes the learner's weighted loss at episode $\ell$.

    By Proposition~\ref{thm:moment-threshold}, when the evaluated model matches the ground truth ($i = i_\star$), $S(y_{(\ell)}; i_\star, j_\ell) = 1$ with probability at least $1 - \delta_1$ whenever $(j_\ell \to i_\star) \in E$, implying $\EE[\hat{Z}_{\ell, i_\star} \mid \mathcal{F}_{\ell-1}] \le \delta_1$.
    A martingale concentration argument bounds $\sum_{\ell=1}^L \hat{Z}_{\ell, i_\star}$, which in turn upper-bounds $\sum_{\ell=1}^L A_\ell$.
    And another application of martingale concentration lower-bounds $\sum_{\ell=1}^L A_\ell$ in terms of the conditional expectations $\mu_\ell := \EE[A_\ell \mid \mathcal{F}_{\ell-1}]$.
    And when $i \neq i_\star$, Proposition~\ref{thm:moment-threshold} ensures $S(\vec{y}^{(\ell)}; i, j_\ell) = 0$ with probability at least $1 - \delta_1$ whenever $(j_\ell \to i) \in E$, which yields
    \begin{align*}
      \mu_\ell
      &\geq (1-\delta_1) \sum_{i\neq i_\star}
                \frac{p_{\ell,i}q_{\ell,i}}{q_{\ell,i}+\beta}
      =(1-\delta_1)\left(1-p_{\ell,i_\star}
          -\beta\sum_{i\neq i_\star}
                            \frac{p_{\ell,i}}{q_{\ell,i}+\beta}\right).
    \end{align*}
    By a graph-theoretic lemma (Lemma~\ref{thm:graph-lemma}), the last sum satisfies $\sum_{i\neq i_\star} \frac{p_{\ell,i}}{q_{\ell,i}+\beta} = O(\bar\alpha \log(1 + \frac{1}{\bar\alpha\beta}))$.
    Combining the upper and lower bounds on $\sum_{\ell=1}^L A_\ell$ leads to an upper bound on the cumulative probability of selecting nonmatching controllers, $\sum_{\ell=1}^L (1 - p_{\ell, i_\star})$.
    Finally, applying Azuma-Hoeffding's inequality shows that $j_\ell = i_\star$ for a strict majority of episodes, guaranteeing successful identification of the matching controller, i.e. $\hat{i} = i_\star$.
\end{proof}

In Appendix~\ref{app:simulation}, we illustrate Algorithm~\ref{alg:meta-graph} on a second-order system with an unknown gain~\citep{anderson2000multiple}.
In this example, the conditions of Definitions~\ref{assum:observability} and~\ref{def:probing} hold for every model-controller pair, so the feedback-information graph is complete and $\alpha = 1$.
As shown in Figure~\ref{fig:simulation}, because all candidates are evaluated simultaneously with provable finite-time bounds, Algorithm~\ref{alg:meta-graph} achieves a much smaller transient than both classical estimator-based supervision and our prior algorithm~\citep{sun2026online}, which evaluates only the model associated with the active controller.

\section{Conclusion}
In this work, we presented a non-asymptotic framework for supervisory switching control of unknown, partially observed linear dynamical systems with a collection of $N$ candidate controllers.
Our approach introduces a cross-model evaluation criterion that combines observability-based state estimation and a short-horizon moment-based threshold test along carefully chosen probing directions.
The criterion enables simultaneous certification or falsification of multiple candidate models even under potentially destabilizing controllers.
By structuring the permissible model-controller pairs as a directed feedback-information graph, our online switching algorithm provably identifies the matching stabilizing controller within $O(\alpha \, \mathrm{polylog}(N))$ steps while preserving finite closed-loop $\mathcal{L}_2$-gain stability.
This result improves sample complexity by up to a factor of $N$ over prior single-model evaluation approaches.

An important future direction is extending these non-asymptotic guarantees to nonlinear dynamical systems.
While our current analysis relies on linear observability matrices and linear input-output Markov parameters, extending our current method to nonlinear plants introduces challenges in analyzing nonlinear observer techniques and local model distinguishability.

\clearpage

\bibliographystyle{plainnat}
\bibliography{reference}

\clearpage

\appendix

\section{Additional related work}
\label{sec:related-work}

\paragraph{Classical methods on supervisory switching control.}
 
Estimator-based supervision~\citep{morse1997supervisory,hespanha2001tutorial,hespanha2003hysteresis,narendra1997adaptive} maintains a bank of candidate models together with their associated controllers, so that a matched pair of model and controllers provides satisfactory performance.
In practice, a finite covering of the plant uncertainty is used to balance approximation accuracy against computational cost~\citep{morse1996supervisory,anderson2000multiple}.
Crucially, the stability argument for this approach only requires the estimation error of the true model to remain small, so the supervisor may settle on a mismatched candidate as long as its model remains consistent with the observed data~\citep{hespanha2001tutorial}.
Consequently, while the multi-estimator evaluates multiple models simultaneously, this approach does not quantify the time needed to find a suitable controller.
Our work can be viewed as a non-asymptotic development of the multi-estimator, where we ensure the performance of cross-model evaluations when a finite trajectory provably carries enough information to support it.
 
A complementary line of work assesses candidate controllers directly through closed-loop performance rather than through model predictions.
In unfalsified control~\citep{safonov1997unfalsified,stefanovic2008safe}, each candidate is evaluated through a virtual reference signal, namely the reference that would have produced the observed data had that candidate been in the loop, so that controllers can be falsified without ever being inserted into the feedback-loop.
Later works designed improved tests that can better detect destabilizing controllers before applying them to the real system~\citep{baldi2010multi,dehghani2009validating}.
Similar to the estimator-based approach, the guarantees of these unfalsified control methods are asymptotic in nature.
 
\paragraph{Learning-based control.}
Learning-based control has received significant attention in the literature over the recent years.
One popular paradigm within this space identifies the unknown dynamics parameters and then applies the certainty equivalence principle.
Recent advances in non-asymptotic system identification have extended classical asymptotic results~\citep{ljung1998system} to quantify the sample complexity of learning linear systems from a single trajectory~\citep{simchowitz2018learning,ziemann2023tutorial,li2024learning}, including partially observed systems whose parameters are reconstructed via the Ho-Kalman algorithm~\citep{oymak2019non,sarkar2021finite,bakshi2023new}.
Together with certainty equivalence, these estimation guarantees yield rigorous performance bounds for robust controller design~\citep{dean2020sample} and efficient regret for linear quadratic control~\citep{mania2019certainty,simchowitz2020naive,tian2023can}.
These methods are designed for settings in which no prior on the true model is known in advance, and their analyses presuppose open-loop stability or a known stabilizing controller.
Our prior work~\citep{sun2024least,sun2026online} built on this machinery, using least-squares estimates to distinguish the candidate controllers, whereas the criterion developed in this paper dispenses with system identification and tests each candidate model directly against the observed trajectory.
 
The second thread treats the control of unknown systems as an online decision-making problem.
The results leveraged online learning techniques such as optimism~\citep{abbasi2011regret}, posterior sampling~\citep{ouyang2017learning}, online convex optimization~\citep{agarwal2019online,chen2021black}, or multiplicative weights update~\citep{li2023online,muehlebach2026the}, and their performance is supported by general complexity measures for interactive decision making~\citep{du2021bilinear,foster2022model}.
However, these analysis still typically impose stability or boundedness assumptions and, in the case of episodic RL~\citep{jin2018q}, require independence of decisions through periodic resets.
 
\paragraph{Bandits with side observations.}
 
Our feedback-information graph connects the supervisory control problem to online learning with side observations~\citep{mannor2011bandits}, where selecting an action also reveals the losses of its out-neighbors in a directed feedback graph.
This framework interpolates between bandit feedback, where each action reveals only its own loss, and full-information feedback, where every action reveals all losses, and its regret is governed by graph-theoretic quantities such as the independence number of the feedback graphs~\citep{alon2015online,alon2017nonstochastic,russo2025pure}.
Our algorithm builds on the importance-weighted multiplicative weights updates developed for this setting, particularly the EXP3-IX algorithm~\citep{kocak2014efficient,lattimore2019bandit}.
 
\section{Mathematical preliminaries}
\label{sec:prelim}
\subsection*{Notations}
For a matrix $M$, we denote $\norm{M}_F$ as its Frobenius norm, $\opnorm{M} = \sigma_{\max}(M)$ as its operator norm, $\rho(M)$ as its spectral radius, and $\tr(M)$ as its trace.
For a stable linear system $(C, A, B)$, we define its $\mathcal{H}$-infinity norm as $\hinfnorm{C, A, B} = \sup_{\norm{s}=1} \sigma_{\max}(C(s I - A)^{-1}B)$, and for simplicity, we use shorthand $\hinfnorm{C, A} = \hinfnorm{C, A, I}$ and $\hinfnorm{A} = \hinfnorm{I, A, I}$.

To streamline our exposition, we occasionally omit constant factors that do not meaningfully impact our conclusions.
We use big-$O$ notation $f \in O(g)$ if $\limsup_{x \to \infty} f(x) / g(x)$ is finite.
Also, we write $f \lesssim g$ if $f \le c \cdot g$ for some universal constant $c$, and $f \asymp g$ if both $f \lesssim g$ and $g \lesssim f$ hold.
Unless stated otherwise, all dimension-dependent terms are written explicitly.
For example, we consider the trace and Frobenius norm to be dimension-dependent, whereas the operator norm is dimension-free.

\subsection*{Technical lemmas}
To analyze stochastic disturbances, we consider a generalization of Gaussian random variables.
Intuitively, a \textit{sub-Gaussian random variable} has tails that are dominated by a Gaussian distribution.
\begin{definition}
    We say that a zero-mean random vector $X \in \RR^d$ is $\sigma^2$-\textit{sub-Gaussian} if for every unit vector $v$ and real value $\lambda$, it holds that
    \[\EE \exp(\lambda \inner{v}{X}) \le \exp(\sigma^2\lambda^2 / 2).\]
\end{definition}
Sub-Gaussian random variables satisfy the following concentration inequality.
\begin{lemma}
\label{thm:sub-gaussian-mgf}
    A real-valued $\sigma^2$-sub-Gaussian random variable $X$ satisfies
    \[\EE[\exp(X^2/(2\sigma)^2)] \le \sqrt{2}.\]
\end{lemma}
\begin{proof}
    Let \(Z\sim \mathcal{N}(0,1)\), and set
\[
Y=\frac{Z}{\sqrt2\,\sigma}.
\]
Then, using the definition of sub-Gaussian r.v.,
\[
\mathbb E_X[e^{YX}\mid Z]\le e^{\sigma^2Y^2/2}
= e^{Z^2/4}.
\]
Also, by the Gaussian moment generating function,
\[
\mathbb E_Z\left[e^{YX} \mid X \right]
=
\mathbb E_Z\left[e^{ZX/(\sqrt2\sigma)}\right]
=
e^{X^2/(4\sigma^2)}.
\]
Therefore,
\[
\begin{aligned}
\mathbb E\left[e^{X^2/(4\sigma^2)}\right]
&=
\mathbb E_X\mathbb E_Z\left[e^{ZX/(\sqrt2\sigma)}\right] \\
&=
\mathbb E_Z\mathbb E_X[e^{YX}] \\
&\le
\mathbb E_Z\left[e^{Z^2/4}\right] \\
&= \frac{1}{\sqrt{2\pi}} \int_\RR \exp(-z^2/4) \, dz 
= \sqrt2.
\end{aligned}
\]
Thus
\[
\mathbb E\left[e^{X^2/(2\sigma)^2}\right]\le \sqrt2.
\]
\end{proof}

Next, we show a graph-theoretic lemma that controls the expected loss incurred in each episode of Algorithm~\ref{alg:meta-graph}.
\begin{lemma}
\label{thm:graph-lemma}
Consider a directed graph $G = ([N], E)$ that includes all self-loops.
Let $\alpha$ be its independence number, which is defined as the cardinality of the largest vertex subset $I \subseteq [N]$ such
that no two distinct vertices in $I$ are connected by a directed edge in $E$, and choose any $\bar\alpha \ge \alpha$.
Let $p$ be any probability distribution on $[N]$, and set
$q_i=\sum_{j:\,j\to i}p_j$. 
For any $\beta > 0$,
\begin{equation}
  \sum_{i=1}^N\frac{p_i}{q_i+\beta}
  \leq 2 \bar\alpha \log\!\left(1+\frac{1}{2 \bar\alpha \beta}\right).
  \label{equ:graph-bound}
\end{equation}
\end{lemma}

\begin{proof}
First fix a subset $S\subseteq[N]$ and put $s=\sum_{i\in S}p_i$.
We claim
\begin{equation}
  \sum_{i\in S}p_i^2
  +2\!\sum_{(j \to i) \in E(S)}p_ip_j
  \geq\frac{s^2}{\bar\alpha}.
  \label{equ:weighted-independence}
\end{equation}
The case $s=0$ is immediate. 
Otherwise, ignore zero-probability vertices and assign an independent exponential r.v. $x_i \sim \mathsf{Exp}[p_i]$ to each remaining vertex $i\in S$.
Let $N_{S, i}$ be the set of all undirected neighbors of $i$ in $S$.
We select a vertex $i$ if $x_i \le x_j$ for every $j \in N_{S, i}$, so that the selected vertices form an independent set.
Writing
\[
  r_i:=p_i+\sum_{j \in N_{S, i}}p_j.
\]
We note that the probability that we select vertex $i$ is equal to
\begin{align*}
    \int_{0}^\infty p_i e^{-p_i x} \prod_{j \in N_{S, i}} e^{-p_j x}\, dx
    = \int_{0}^\infty p_i e^{-r_i x} \, dx
    = \frac{p_i}{r_i}.
\end{align*}
Hence, $\sum_{i\in S} p_i/r_i \leq \alpha \le \bar\alpha$.
By Cauchy--Schwarz,
\[
  s^2
  \leq
  \left(\sum_{i\in S}p_i r_i\right)
  \left(\sum_{i\in S}\frac{p_i}{r_i}\right)
  \leq \bar\alpha \sum_{i\in S}p_i r_i,
\]
which proves~\eqref{equ:weighted-independence}.

Therefore,
\begin{equation}
\begin{aligned}
  \sum_{i\in S}p_iq_i
  &= \sum_{i \in S} p_i^2 + \sum_{i \in S} \sum_{j : (j \to i) \in E} p_i p_j \\
  &\ge \sum_{i\in S}p_i^2 +\sum_{(j \to i) \in E(S)} p_i p_j
   \ge\frac{s^2}{2\bar\alpha}.
\end{aligned}
\label{equ:directed-subset}
\end{equation}
For $u\geq0$, let $F(u)=\sum_{i:q_i\leq u}p_i$.
Applying \eqref{equ:directed-subset} to $S=\{i:q_i\leq u\}$ gives
\[
  \frac{F(u)^2}{2\bar\alpha}\leq uF(u) \implies F(u)\leq\min\{1,2\bar\alpha u\}.
\]
It follows that
\begin{align*}
  \sum_{i=1}^N\frac{p_i}{q_i+\beta}
  &=\sum_i p_i\int_{q_i}^{\infty}\frac{du}{(u+\beta)^2}
    =\int_0^\infty\frac{F(u)}{(u+\beta)^2}\,du\\
  &\le \int_0^\infty \frac{\min\{1,2\bar\alpha u\}}{(u+\beta)^2}\,du \\
  &\le \int_0^{1/(2\bar\alpha)} \frac{2\bar\alpha u}{(u+\beta)^2}\,du + \int_{1/(2\bar\alpha)}^\infty \frac{1}{(u+\beta)^2}\,du \\
  &= 2\bar\alpha \bigg[\frac{\beta}{u + \beta} + \log(\beta + u) \bigg]_{u=0}^{1/(2\bar\alpha)} + \frac{1}{1/(2\bar\alpha) + \beta}\\
  &=2\bar\alpha\log\left(1+\frac1{2\bar\alpha\beta}\right).
\end{align*}
\end{proof}
 
\section{Numerical simulation}\label{app:simulation}

To illustrate the behavior of Algorithm~\ref{alg:meta-graph}, we consider an example studied in \citet{anderson2000multiple}, which we also used in our prior work \citep{sun2026online}.
In this example, we have a second-order system with an unknown gain $k$,
\[ \left\{\frac{k(s-1)}{(s+1)(s-2)} : k \in [1, 40)\right\}. \]
This family of systems is difficult to control because every member of this family is non-minimum phase and because the uncertainty on $k$ spans a large range.
We construct candidate models on the finite covering $k_i = 1.2^i$, $i \in \{0, \dots, 20\}$, which is found to closely approximate the uncertainty on $k$ \citep{anderson2000multiple}. This gives $N = 21$ candidate models and controllers.

We simulate a zero-order-hold discretization of this system at 1000 Hz.
To illustrate a \textit{worst-case} scenario, we sample a large initial condition sampled from $\mathcal{N}(0, 20)$, set the true index to be $i_\star = 20$, and fix every policy's starting index to be $i_1 = 0$.
We compare three approaches:
\begin{itemize}
    \item \textit{Classical estimator-based supervision}, with a dwell time of 0.2 seconds and a forgetting factor of $0.04$.
    \item \textit{Online algorithm with single-model evaluation} \citep{sun2026online}, which in each episode evaluates only the model associated with the active controller. We set the maximum episode length to 0.2 seconds, use an early-stopping rule when the active controller is determined to be destabilizing, and use an exploration factor of $0.3\ell$ for the $\ell$th episode.
    \item \textit{Graph-structured online algorithm with cross-model evaluations (Algorithm~\ref{alg:meta-graph})}, with a maximum episode length of 0.2 seconds and parameters $\beta = 1/80$, $\eta = 1.0$. We also add an early-stopping rule to the evaluation criterion (Algorithm~\ref{alg:moment-threshold}), so that we switch as soon as the rejection threshold associated with the active controller is exceeded.
\end{itemize}
Since our theoretical bounds are quite conservative, all policy parameters were tuned to give the best empirical behavior for this scenario.

For this example, the feedback-information graph that meets Definition~\ref{def:graph} is complete.
The probing horizon is only $h = 10$ steps of 1 ms each.
Over such a short window, even a destabilizing closed loop cannot amplify the residual by much.
Therefore, even though every mismatched model-controller results in an unstable closed-loop system, the conditions in Definition~\ref{def:probing} still holds for every model-controller pair.
At the same time, the distinct gains $k_i$ keep the Markov parameters of the possible closed-loop models separated.
Consequently, the independence number is $\alpha = 1$, and every episode certifies or falsifies all $21$ candidate models at once, regardless of which controller is active.

\begin{figure}[!t]
    \flushright
    \includegraphics[width=0.97\linewidth]{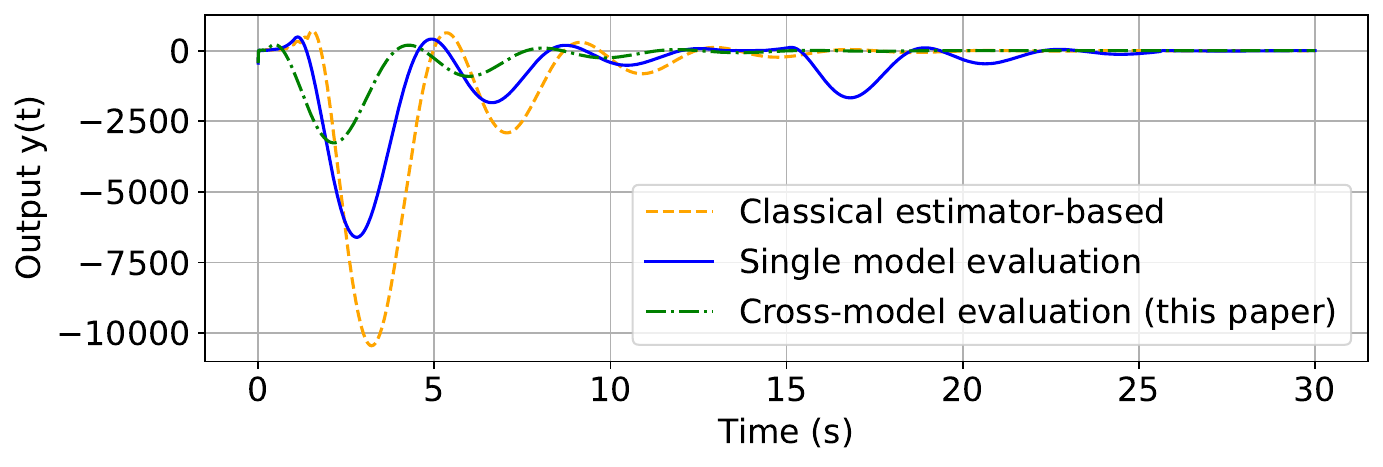}\\[0.2em]
    \centering
    {\small (a) Plant output $y(t)$.}\\[1.0em]
    \flushright
    \includegraphics[width=0.90\linewidth]{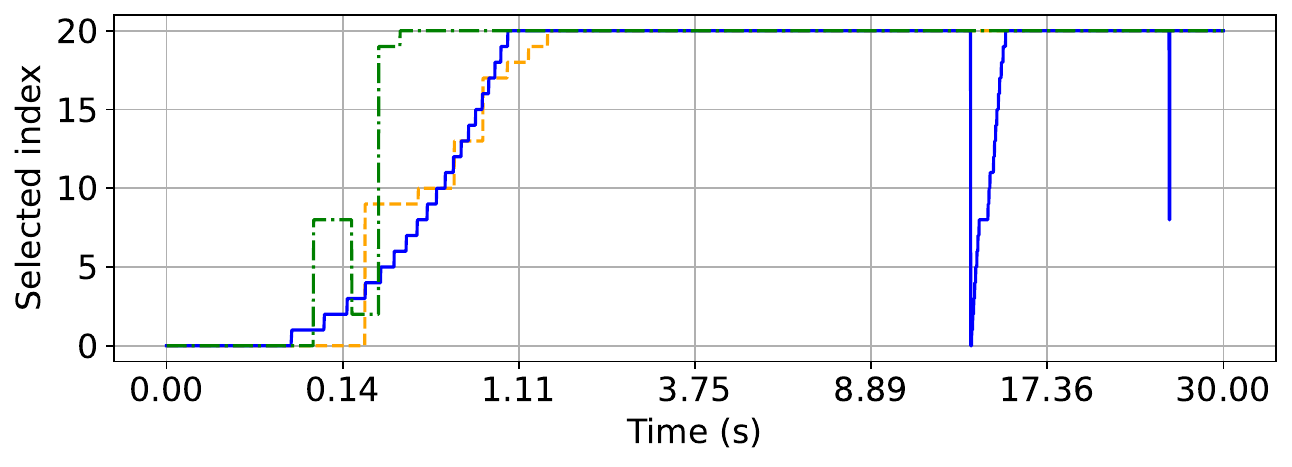}\\[0.2em]
    \centering
    {\small (b) Index of the candidate controller selected by each policy (true index $i_\star = 20$).}
    \caption{Numerical simulation comparing Algorithm~\ref{alg:meta-graph} (green, dash-dotted) against classical estimator-based supervision (orange, dashed) and the online algorithm with single-model evaluation of \citet{sun2026online} (blue, solid). The time axis in (b) is cube-root scaled to clearly show every switching decision.}
    \label{fig:simulation}
\end{figure}

Figure~\ref{fig:simulation} shows the resulting trajectories from applying each of the three approaches.
Algorithm~\ref{alg:meta-graph} produces a much smaller initial transient: its output peaks at $|y| \approx 3.3 \times 10^3$, compared with $6.6 \times 10^3$ for our prior online algorithm and $1.0 \times 10^4$ for the estimator-based supervisor.
Two features of our method account for this gap.
First, because all candidates are evaluated simultaneously, the weights concentrate on the matching controller within a few episodes. 
In contrast, the online algorithm with single evaluation learns about only one model per episode and has to deploy all $21$ candidates in turn.
Second, our evaluation criterion has provable guarantees that can support the usage of explicit stopping rules, but the estimator-based supervisor is limited by a dwell time that has no precise theoretical backing.
For these two reasons, Algorithm~\ref{alg:meta-graph} reaches $i_\star$ more quickly than the other two approaches.
Since the state can grow quickly while a mismatched controller is in the loop, the shorter identification phase of Algorithm~\ref{alg:meta-graph} translates directly into a markedly smaller peak.

The advantage of our cross-model evaluation can be seen beyond the initial peak.
In particular, the algorithm in \citet{sun2026online} can only evaluate a candidate by applying its associated controller to the loop, so it requires an exploration factor that periodically forces it to revisit the other candidates.
In Figure~\ref{fig:simulation}(b), this explicit exploration happens around $t = 15$ seconds, when the algorithm sweeps through all candidates once more and causes a second transient.
In comparison, Algorithm~\ref{alg:meta-graph} requires no such explicit exploration.
With a complete feedback-information graph, every candidate is certified or falsified for every episode, so exploration comes at no extra cost and the algorithm does not need to deviate from $i_\star$ after the initial exploration phase.

\section{Proof of Proposition~\ref{thm:moment-threshold}}
\label{sec:proof-criteria}
We shall show that the conditions in Proposition~\ref{thm:moment-threshold} are met with
\begin{align*}
    \tau_2 = 32h \left(\frac{\Lambda^2(\sigma_r^2 + \sigma_w^2) + \sigma_\eta^2}{\sigma_u^2 \gamma^2} \log(2\sqrt{2} N /\delta_1) + \frac{1}{4}\right)\log(2N/\delta_1).
\end{align*}
For convenience, let $\BB$ be the $\sigma$-algebra on all process and observation noise $w$'s and $\eta$'s and $\mathcal{G}_t$ be the $\sigma$-algebra on $\BB$ and all exploratory signal $u$'s up to time $t-1$. 
Define interlacing classes $\mathcal{C}_c, c = 0, \dots, h-1$, such that
\[\mathcal{C}_c = \{h + \tau_1 < t \le h + \tau_1 + \tau_2 \mid t \equiv c \hspace{-0.5em} \pmod{h}\}.\]

\subsection{Bounding state estimation error}
We first show that the state estimation procedure described in Section~\ref{sec:state-estimation} and lines 1-5 in Algorithm~\ref{alg:moment-threshold} produces an accurate estimate $\hat{x}_{t-h}$ for each time step $t$.

\begin{lemma}
    Consider an evaluated model $i$ that matches the true system ($i = i_\star$) and active controller $j$, and $M_{[i,j]}$ is $\varepsilon_c$-strictly observable.
    Define constants $\kappa_1 > 0, \kappa_2 > 0, \lambda_1 \ge 1, \lambda_2 \in (0, 1)$, and $\Theta > 0$ so that
    \begin{align*}
        &\opnorm{C_{[i, j]}} < \kappa_1, \quad \opnorm{A_{[i, j]}} \le \lambda_1, \\
        &\opnorm{(A_{[i, j]} - L_{[i, j]}C_{[i, j]})^s} \le \kappa_2 \lambda_2^s \;\; \forall s, \quad \hinfnorm{A_{[i, j]} - L_{[i, j]}C_{[i, j]}} \le \sqrt{\Theta}.
    \end{align*}
    Then, for any time step $s$ and 
    \[ \tau_1 = \max\left\{\nu, \frac{\log(\varepsilon_c^{-1}\kappa_1\kappa_2\nu\lambda_1^{\nu-1})}{\log(1/\lambda_2)} \right\} \in O(\nu),\]
    the state estimate $\hat{x}_{s + \tau_1}$ computed from the observations $y_{s}, \dots, y_{s + \tau_1 - 1}$ satisfies that the difference
    $x_{s+\tau_1} - \hat{x}_{s + \tau_1}$ is $((1 + \Theta)\sigma_w^2 + (1 + \Theta\opnorm{L_{[i,j]}}^2)\sigma_\eta^2)$-sub-Gaussian and independent of the exploratory signal $u$.
\end{lemma}

\begin{proof}
    From \eqref{equ:obs-dynamics-alg}, the initial estimation error is equal to
    \[x_s - \hat{x}_s = \mathcal{O}_{[i,j]}^\dagger \left(
    \begin{bmatrix} 
    C_{[i,j]} & C_{[i,j]}A_{[i,j]} & C_{[i,j]}A^2_{[i,j]} & \dots & C_{[i,j]}A^{\nu-2}_{[i,j]} \\
    0 & C_{[i,j]} & C_{[i,j]}A_{[i,j]} & \dots & C_{[i,j]}A_{[i,j]}^{\nu-3} \\
    0 & 0 & C_{[i,j]} & \dots & C_{[i,j]}A_{[i,j]}^{\nu-4}\\
    \vdots & \vdots & \vdots & \ddots & \vdots \\
    0 & 0 & 0 & \dots & C_{[i,j]} \\
    0 & 0 & 0 & \dots & 0
\end{bmatrix}
    \begin{bmatrix} w_{\nu-1} \\ \vdots \\ w_2 \\ w_1  \end{bmatrix}
    +
    \begin{bmatrix} \eta_\nu \\ \vdots \\ \eta_2 \\ \eta_1  \end{bmatrix}
    \right).\]

We note that the operator norm of the row $[C_{[i,j]} \; C_{[i,j]}A_{[i,j]} \; \dots \; C_{[i,j]}A^{k-1}_{[i,j]}]$ is at most $\kappa_1 k \lambda_1^{k-1}$.
And under observability, $\opnorm{\mathcal{O}_{[i,j]}^\dagger} \le 1/\varepsilon_c$.
This means that $x_s - \hat{x}_s$ is $((\varepsilon_c^{-1} \kappa_1 \nu \lambda_1^{\nu-1})^2(\sigma_w^2+\sigma_\eta^2))$-sub-Gaussian.

With the state observer $L_{[i,j]}$, the difference after $\tau_1$ steps becomes:
\[
x_{s+\tau_1} - \hat{x}_{s+\tau_1} = (A_{[i,j]} - L_{[i,j]} C_{[i,j]})^{\tau_1} (x_s - \hat{x}_s) + \sum_{k=1}^{s+\tau_1-1} (A_{[i,j]} - L_{[i,j]} C_{[i,j]})^{t-k-1} (w_k - L_{[i,j]} \eta_k).
\]
We note that all of the terms on RHS are free of the exploratory signal $u$.
The first term is then $((\varepsilon_c^{-1} \kappa_1\kappa_2 \nu \lambda_1^{\nu-1} \lambda_2^{\tau_1})^2(\sigma_w^2 + \sigma_\eta^2))$-sub-Gaussian.
With our choice of $\tau_1$, this becomes $(\sigma_w^2 + \sigma_\eta^2)$-sub-Gaussian.
And the second term is $(\Theta(\sigma_w^2 + \opnorm{L_{[i,j]}}^2\sigma_\eta^2))$-sub-Gaussian because any finite-dimensional Toeplitz matrix generated by $A_{[i,j]} - L_{[i,j]} C_{[i,j]}$ has operator norm at most $\hinfnorm{A_{[i,j]} - L_{[i,j]} C_{[i,j]}}$ (see \citet{tilli1998singular}).
Summing up these two parts yields the desired claim.
\end{proof}

Therefore, for each time step $t$, we have a state estimate $\hat{x}_{t-h}$ whose error $x_{t-h} - \hat{x}_{t-h}$ is $\sigma_r^2$-sub-Gaussian, where $\sigma_r^2 = (1 + \Theta)\sigma_w^2 + (1 + \Theta\opnorm{L_{[i,j]}}^2)\sigma_\eta^2$.

\subsection{Bounding error probability in the case $i = i_\star$}
In this part, we bound the probability that Algorithm~\ref{alg:moment-threshold} incorrectly rejects the evaluated model $i$ when $i = i_\star$.

\paragraph{Step 1:}
From line 3 of Algorithm~\ref{alg:moment-threshold}, the residual $\widetilde{y}_t$ can be expressed as
\[\widetilde{y}_t = C_{[i,j]} A_{[i,j]}^h (x_{t-h} - \hat{x}_{t-h}) + \sum_{s=t-h}^{t-1} C_{[i,j]} A_{[i,j]}^{t-1-s} w_s + \eta_t.\]
From the previous part, we know that $x_{t-h} - \hat{x}_{t-h}$ is $\sigma_r^2$-sub-Gaussian.
Under Definition~\ref{def:probing}, the quantity $r_{i', t} = a_{i'}^\top \widetilde{y}_t$ is $\sigma^2$-sub-Gaussian for $\sigma^2 = \Lambda^2(\sigma_r^2 + \sigma_w^2) + \sigma_\eta^2$.

We claim that with $M = 4\sigma^2 \log(2\sqrt{2} N / \delta_1)$,
\[ \PP\left(\frac{1}{\tau_2} \sum_{t=h+\tau_1+1}^{\tau} r_{i',t}^2 > M \right) \le \frac{\delta_1}{2N}.\]
From Jensen's inequality, we have
\[\exp\left(\frac{1}{4\sigma^2\tau_2}\sum_t r_{i',t}^2\right)\le \frac{1}{\tau_2} \sum_t \exp(r_{i',t}^2/(2\sigma)^2).\]
Taking expectations yields
\[
 \EE\left[\exp\left(\frac1{4\sigma^2\tau_2}\sum_t r_{i',t}^2\right)\right]
 \le\frac{1}{\tau_2} \sum_t\EE\left[\exp(r_{i',t}^2/(2\sigma)^2)\right] \le \sqrt2,
\]
where the last inequality follows from Lemma~\ref{thm:sub-gaussian-mgf}.
Then, from Markov's inequality, we have
\[ \PP\left(\frac{1}{\tau_2} \sum_{t=h+\tau_1+1}^{\tau} r_{i',t}^2 > M \right) \le \sqrt2 e^{-M/(4\sigma^2)}.\]

With our choice on $M$,
\begin{equation}
\label{equ:energy-error}
 \PP(\Phi_{i'}^{(1)} > M) \le \sqrt2 e^{-M/(4\sigma^2)} = \delta_1/(2N).
\end{equation}
Then, by union bound, the probability that $\Phi_{i'}^{(1)} \le M$ for all $i'$ is at least $1 - \delta_1/2$.

\paragraph{Step 2:}
Conditioned on $\BB$, all $r_{i', t}$'s are fixed.
Within an interlacing class $\mathcal C_c$, the sliding window $z_t$'s for $t \in \mathcal C_c$ are disjoint.
Therefore, for every real $\lambda$, applying the tower property on $\mathcal{G}_t, t \in \mathcal{C}_c$, yields
\begin{align*}
 \EE\left[\exp\left(\lambda\sum_{t\in\mathcal C_c}r_{i', t}\bar{z}_{i',t}\right) \mid\BB\right]
 \le \exp\left(\frac{\lambda^2\sigma_u^2}{2} \sum_{t\in\mathcal C_c}r_{i', t}^2\right).
\end{align*}
Applying Conditional H\"older's inequality across the $h$ classes yields
\begin{equation}
\begin{aligned}
 \EE\left[\exp\left(\lambda\sum_{t}r_{i', t}\bar{z}_{i',t}\right) \mid\BB\right]
 &\le \prod_{c=0}^{h-1} \EE\left[\exp\left(h\lambda\sum_{t\in\mathcal C_c}r_{i', t}\bar{z}_{i',t}\right) \mid\BB\right]^{1/h} \\
 &\le\exp\left(\frac{h\lambda^2\sigma_u^2}{2} \sum_{t}r_{i', t}^2\right).
\end{aligned}
 \label{eq:target-product-mgf}
\end{equation}

\paragraph{Step 3:}
On the event $\mathcal{E}$ that $\Phi_{i'}^{(1)} \le M$ for all $i'$, \eqref{eq:target-product-mgf} is at most
$\exp(h\lambda^2 \sigma_u^2 \tau_2 M/2)$.
Applying Markov's inequality and setting $\lambda=\gamma/(2h M)$ yields
\begin{equation}
\begin{aligned}
 \PP\left(\mathcal{E} \wedge \frac{1}{\tau_2} \sum_t r_{i', t}\bar{z}_{i', t} > \frac{\sigma_u^2 \gamma}{2} \right)
 &\le \inf_{\lambda> 0} \exp\left(h\lambda^2 \sigma_u^2 \tau_2 M/2 - \lambda \sigma_u^2 \tau_2 \gamma / 2\right) \\
 &\le\exp\left(-\frac{\tau_2 \sigma_u^2 \gamma^2}{8 h M}\right).
\end{aligned}
 \label{eq:target-product-tail}
\end{equation}
By choosing
\[ \tau_2 \ge \frac{32h (\Lambda^2(\sigma_r^2 + \sigma_w^2) + \sigma_\eta^2)}{\sigma_u^2\gamma^2} \log(2\sqrt{2} N/\delta_1)\log(2N/\delta_1) = \frac{8hM}{\sigma_u^2\gamma^2} \log(2N/\delta_1) ,\] the final term in \eqref{eq:target-product-tail} becomes $\delta_1/(2N)$.
A union bound on \eqref{equ:energy-error} and \eqref{eq:target-product-tail} establishes that when $i = i_\star$, Algorithm~\ref{alg:moment-threshold} correctly certifies model $i$ with probability at least $1-\delta_1$.

\subsection{Bounding error probability in the case $i \neq i_\star$}
In this part, we bound the probability that Algorithm~\ref{alg:moment-threshold} incorrectly accepts the evaluated model $i$ when $i \neq i_\star$.
In this case, $r_{i_\star, t} = a_{i_\star}^\top \Delta_{i_\star, i} z_t + v_t$, where $a_{i_\star}^\top \Delta_{i_\star, i} = \norm{a_{i_\star}^\top \Delta_{i_\star, i}} b_{i_\star}^\top$ and $v_t$ are noise independent of $z_t$.
So, the correlation between $r_{i_\star, t}$ and $\bar{z}_{i_\star, t}$ is equal to $\rho := \norm{a_{i_\star}^\top \Delta_{i_\star, i}} > \gamma$.

\paragraph{Step 1:}
For any fixed $\lambda \ge 0$, directly integrating over the Gaussian density on $\bar{z}_{i_\star, t}$ gives
\begin{equation}
\label{eq:rival-one-step}
\begin{aligned}
 \EE \left[\exp\left(-\frac{\lambda r_{i_\star, t} \bar{z}_{i_\star, t}}{\sigma_u^2} -\frac{\lambda^2 r_{i_\star, t}^2}{2\sigma_u^2}\right) \mid \mathcal{G}_{t-h}\right]
 &=\frac1{\sqrt{2\pi}\sigma_u} \int_{\RR} \exp\left(-\frac{(\bar{z}+\lambda(\rho \bar{z}+v))^2}{2\sigma_u^2}\right)\,d\bar{z}\\
 &=\frac1{\sqrt{2\pi}\sigma_u} \int_{\RR} \exp\left(-\frac{((1+\lambda \rho)\bar{z} + \lambda v)^2}{2\sigma_u^2}\right)\,d\bar{z} \\
 &=\frac{1}{1+\lambda\rho}.
\end{aligned}
\end{equation}

Then, for each interlacing class define
\begin{equation}
 Q_c:=\prod_{t\in\mathcal C_c}
       (1+\lambda\rho) \exp\left(-\frac{\lambda r_{i_\star, t} \bar{z}_{i_\star, t}}{\sigma_u^2} -\frac{\lambda^2 r_{i_\star, t}^2}{2\sigma_u^2}\right).
 \label{eq:class-Q}
\end{equation}
Within an interlacing class $\mathcal C_c$, the sliding window $z_t$'s for $t \in \mathcal C_c$ are disjoint.
So, we can apply tower law and \eqref{eq:rival-one-step} along each class.
It follows that $\EE[Q_c\mid\BB]=1$.
Let
\begin{equation*}
 Q:=\prod_cQ_c
 =(1+\lambda\rho)^{\tau_2} \exp\left(-\sum_t \frac{\lambda r_{i_\star, t} \bar{z}_{i_\star, t}}{\sigma_u^2} - \sum_t \frac{\lambda^2 r_{i_\star, t}^2}{2\sigma_u^2}\right).
 \label{eq:Q}
\end{equation*}
Applying H\"older's inequality gives
\begin{equation}
 \EE[Q^{1/h}\mid\BB]
 \le \prod_c(\EE[Q_c \mid \BB])^{1/h}=1.
 \label{eq:Q-holder}
\end{equation}

\paragraph{Step 2:}
Define the event
\[
 \mathcal{A}:=\left\{\sum_t r_{i_\star, t}^2\le \tau_2 M \quad \wedge \quad \sum_t r_{i_\star, t} \bar{z}_{i_\star, t} \le \tau_2 \frac{\rho\sigma_u^2}{2}\right\}.
\]
Because $\rho > \gamma$, acceptance of model $i$ under truth $i_\star \neq i$ implies $\mathcal{A}$, under which
\begin{equation}
 Q^{1/h} 
 \ge \exp\left(\frac{\tau_2}{h} \left\{\log(1+\lambda\rho) - \frac{\lambda\rho}{2} - \frac{\lambda^2M}{2\sigma_u^2}\right\}\right)
 \ge\exp\left(\frac{\tau_2}{h} \left\{\frac{\lambda\rho}{2} - \frac{\lambda^2(M/\sigma_u^2+\rho^2)}{2}\right\}\right)
 \label{eq:Q-on-acceptance}
\end{equation}
where the last inequality follows from the fact that $\log(1+x)\ge x-x^2/2 \;\; \forall x \ge 0$.
Choosing
\[
 \lambda=\frac{\rho}{2(M/\sigma_u^2 + \rho^2)},
\]
\eqref{eq:Q-holder} and \eqref{eq:Q-on-acceptance} together implies
\begin{equation*}
 \PP(\mathcal{A}\mid\BB)
 \le\exp\left(-\frac{\tau_2 \rho^2}{8h (M / \sigma_u^2 +\rho^2)}\right).
 \label{eq:rival-product-tail}
\end{equation*}
Then, with
\begin{align*}
\tau_2 
&\ge 32h \left(\frac{\Lambda^2(\sigma_r^2 + \sigma_w^2) + \sigma_\eta^2}{\sigma_u^2 \gamma^2} \log(2\sqrt{2} N /\delta_1) + \frac{1}{4}\right)\log(1/\delta_1) \\
&= \left(\frac{8h M}{\sigma_u^2 \gamma^2} + 8h\right) \log(1/\delta_1)
\ge \left(\frac{8h M}{\sigma_u^2 \rho^2} + 8h\right) \log(1/\delta_1),
\end{align*}
the probability of incorrectly accepting model $i$ when $i \neq i_\star$ is at most $\delta_1$.
 
\section{Proof of Theorem~\ref{thm:graph-exp3}}
\label{sec:proof-graph-bandit}
We shall prove the two parts of this theorem in order.

\subsection*{Part 1: Sample complexity bound}

Denote
\[ A_\ell:=\sum_{i=1}^N p_{\ell,i}\widehat Z_{\ell,i}. \]
Recall that, under Proposition~\ref{thm:moment-threshold}, for every $i$ where $j_\ell \to i \in E$,
\begin{align*}
    \PP(S(\vec{y}^{(\ell)}; i, j_\ell) = 1) &\ge 1-\delta_1, \quad \text{if } i = i_\star; \\
    \PP(S(\vec{y}^{(\ell)}; i, j_\ell) = 0) &\ge 1-\delta_1, \quad \text{if } i \neq i_\star.
\end{align*}
For convenience, let us choose
\begin{equation*}
  \delta_1 = \frac{1}{17},
  \qquad
  \beta:=\frac{1}{80\bar\alpha},
  \qquad \eta:=\frac{\beta}{2},
  \qquad
  L = \frac{30}{\beta}\log\frac{3N}{\delta_2},
\end{equation*}
and show that with these constants, Algorithm~\ref{alg:meta-graph} commits to controller $i_\star$ with probability at least $1-\delta_2$.
We note that the numerical constants are conservative and chosen to simplify the proof.

\paragraph{Step 1:}
Let $W_\ell=\sum_i w_{\ell,i}$.
Then,
\[\frac{W_{\ell+1}}{W_\ell} = \frac{\sum_i w_{\ell_i} \exp(-\eta \hat{Z}_{\ell,i})}{W_\ell}\]
Using $e^{-x}\leq1-x+x^2/2$ for $x\geq0$
and $\log(1+x)\leq x$, we obtain
\begin{align*}
  \log\frac{W_{\ell+1}}{W_\ell}
  &= \log \frac{\sum_i w_{\ell_i} \exp(-\eta \hat{Z}_{\ell,i})}{W_\ell} \\
  &\le \log \left(1 + \frac{\sum_i  -\eta w_{\ell,i}\hat{Z}_{\ell,i} + \eta^2 w_{\ell,i} \hat{Z}_{\ell,i}^2 / 2}{W_\ell}\right) \\
  &= \log \left(1 + \sum_i  -\eta p_{\ell_i}\hat{Z}_{\ell,i} + \frac{\eta^2}{2} p_{\ell_i} \hat{Z}_{\ell,i}^2\right) \\
  &\leq -\eta A_\ell + \frac{\eta^2}{2}\sum_i p_{\ell,i}\widehat Z_{\ell,i}^2.
\end{align*}
Also, $W_1=N$ and
\[
  W_{L+1} \geq w_{L+1,i^\star} = \exp{\left(-\eta\sum_{\ell=1}^L \widehat Z_{\ell,i^\star}\right)}.
\]
Combining these inequalities with telescoping sum yields
\[
-\eta\sum_{\ell=1}^L \widehat Z_{\ell,i^\star} - \log N \le \log\frac{W_{L+1}}{W_1} 
\le  \sum_{\ell=1}^L \left(-\eta A_\ell + \frac{\eta^2}{2}\sum_i p_{\ell,i}\widehat Z_{\ell,i}^2\right).
\]
It follows that
\[
  \sum_{\ell=1}^L A_\ell
  \leq \sum_{\ell=1}^L\widehat Z_{\ell,i^\star}
      +\frac{\log N}{\eta}
      +\frac\eta2\sum_{\ell=1}^L\sum_i
                           p_{\ell,i}\widehat Z_{\ell,i}^2.
\]
Since $0\leq\widehat Z_{\ell,i}\leq1/\beta$,
\[
  \sum_i p_{\ell,i}\widehat Z_{\ell,i}^2
  \leq\frac{A_\ell}{\beta}.
\]
Setting $\eta=\beta/2$, we conclude that
\begin{equation}
  \frac{3}{4}\sum_{\ell=1}^L A_\ell
  \leq \frac{2\log N}{\beta}
       +\sum_{\ell=1}^L\widehat Z_{\ell,i^\star}.
  \label{equ:potential}
\end{equation}

\paragraph{Step 2:}
Let $\mathcal F_{\ell-1}$ be the $\sigma$-algebra on all history before $\ell$th episode, and denote 
\[ X_\ell = (1 - S(\vec{y}^{(\ell)}; i_\star, j_\ell)) \; \mathbb I\{j_\ell \to i_\star \in E\}. \]
Convexity implies that $e^{\beta x}\leq 1 + x(e^\beta - 1)$ for $x\in[0,1]$.
Then by Proposition~\ref{thm:moment-threshold},
\begin{align*}
\EE \left[\exp(\beta\widehat Z_{\ell,i^\star}) \mid\mathcal F_{\ell-1} \right] 
&=\EE \left[\exp\left(\frac{\beta X_\ell}{q_{\ell, i_\star} + \beta}\right) \mid\mathcal F_{\ell-1} \right] \\
&\le \EE \left[1 + X_\ell (e^{\beta/(q_{\ell,i_\star} + \beta)} - 1) \mid\mathcal F_{\ell-1} \right] \\
&\le 1 + \delta_1 q_{\ell, i_\star} (e^{\beta/(q_{\ell,i_\star} + \beta)} - 1).
\end{align*}
Moreover, because $\log(1+x)\geq x/(1+x)$ for $x\geq0$, we have $e^{\beta/(q_{\ell, i_\star} + \beta)} \leq \frac{\beta}{q_{\ell, i_\star}} + 1.$
Consequently,
\begin{align*}
  \EE \left[ \exp(\beta\widehat Z_{\ell,i^\star}) \mid\mathcal F_{\ell-1} \right] \leq 1 + \delta_1\beta \leq \exp(\delta_1 \beta).
\end{align*}
Iterating this conditional moment bound and applying Markov's inequality
shows that
\[
  \PP \left(\sum_{\ell=1}^L\widehat Z_{\ell,i^\star} < \delta_1 L + \frac{\log(3/\delta_2)}{\beta} \right) \ge 1 - \frac{\delta_2}{3}.
\]
Under this event, setting $\delta_1 = 1/17$ and applying~\eqref{equ:potential} yields
\begin{equation}
  \sum_{\ell=1}^L A_\ell
  \leq\frac{4}{3 \cdot 17}L + \frac{4}{3\beta}(2\log N + \log(3/\delta_2)).
  \label{equ:A-bound}
\end{equation}

\paragraph{Step 3:}
Let $\mu_\ell=\mathbb E[A_\ell\mid\mathcal F_{\ell-1}]$, and note that
\[
  0 \leq A_\ell \leq\sum_i\frac{p_{\ell,i}}{q_{\ell,i}+\beta}.
\]
By Lemma~\ref{thm:graph-lemma}, the RHS is upper bounded by the quantity $D := 2\bar\alpha \log(1 + \frac{1}{2\bar\alpha\beta})$.
Because $e^{-x} \leq 1-x/2$ for $x\in[0,1]$,
\[
  \mathbb E\!\left[e^{-A_\ell/D}\mid\mathcal F_{\ell-1}\right]
  \leq1-\frac{\mu_\ell}{2D}
  \leq\exp\!\left(-\frac{\mu_\ell}{2D}\right).
\]
Therefore,
\[
  \exp\!\left(
     \sum_{\ell=1}^m\frac{\mu_\ell}{2D}
     -\sum_{\ell=1}^m\frac{A_\ell}{D}
  \right),\qquad m=0,\ldots,L,
\]
form a nonnegative supermartingale starting at one.
Ville's inequality implies that, with probability at least $1 - \delta_2 / 3$,
\begin{equation}
  \sum_{\ell=1}^L\mu_\ell - 2 \sum_{\ell=1}^L A_\ell \le 2 D \log(3/\delta_2).
  \label{equ:mean-concentration}
\end{equation}

\paragraph{Step 4:}
For every $i\neq i^\star$, Proposition~\ref{thm:moment-threshold} implies that
observing controllers gives
\begin{align*}
  \EE[\widehat Z_{\ell,i}\mid\mathcal F_{\ell-1}]
  &=\frac{1}{q_{\ell,i}+\beta} \sum_{j:\,j\to i}p_{\ell,j} \EE [1 - S(\vec{y}^{(\ell)}; i, j) \mid\mathcal F_{\ell-1}]\\
  &\ge (1-\delta_1) \frac{q_{\ell,i}}{q_{\ell,i}+\beta}.
\end{align*}
Then
\begin{align*}
  \mu_\ell
  &= \sum_{i\neq i^\star} p_{\ell, i} \EE[\widehat Z_{\ell,i}\mid\mathcal F_{\ell-1}]\\
  &\geq (1-\delta_1) \sum_{i\neq i^\star}
            \frac{p_{\ell,i}q_{\ell,i}}{q_{\ell,i}+\beta}\\
  &= (1-\delta_1)\left(1-p_{\ell,i^\star} -\beta\sum_{i\neq i^\star} \frac{p_{\ell,i}}{q_{\ell,i}+\beta}\right)\\
  &\ge (1-\delta_1)(1-p_{\ell,i^\star}-\beta D).
\end{align*}
Define the cumulative conditional mass on incorrect controllers by
\[
  P_L:=\sum_{\ell=1}^L(1-p_{\ell,i^\star}).
\]
Then, with $\delta_1 = 1/17$,
\begin{equation}
  P_L \le \frac{17}{16} \sum_{\ell=1}^L\mu_\ell+\beta DL.
  \label{equ:mass-from-mean}
\end{equation}
Then,
\eqref{equ:A-bound}, \eqref{equ:mean-concentration}, and \eqref{equ:mass-from-mean} imply that with probability $1 - 2\delta_2/3$,
\begin{equation}
\begin{aligned}
    P_L
    &\le \frac{17}{8} \sum_{\ell=1}^L A_\ell + \frac{17}{8} \frac{\beta D \log(3/\delta_2)}{\beta} + \beta D L\\
    &\le \frac{L}{6} + \frac{17}{6\beta} (2\log N + \log(3/\delta_2)) + \frac{17}{8} \frac{\beta D \log(3/\delta_2)}{\beta} + \beta D L\\
    &\le \left(\frac{1}{6} + \beta D\right) L + \frac{\log(3N/\delta_2)}{\beta} \left(\frac{17}{3} + \frac{17 \beta D}{8}\right).
\end{aligned}
  \label{equ:mass-explicit}
\end{equation}
With the choice $\beta = \frac{1}{80\bar\alpha}$, it follows that
\[\beta D \le 2 \bar\alpha \beta \log\left(1+\frac{1}{2\bar\alpha\beta}\right) = \frac{\log(41)}{40} < 1/10.\]
Then, with $L = 30 \log(3N/\delta_2) / \beta$, we conclude that
\begin{equation}
  \frac{P_L}{L}
  \leq \frac{1}{6} + \frac{1}{10} + \frac{1}{30} \left(\frac{17}{3} + \frac{17}{80}\right) < 0.47.
  \label{equ:mass-bound}
\end{equation}

\paragraph{Step 5:}
Let
\[
  B_L:=\sum_{\ell=1}^L \mathbb{I} \{j_\ell\neq i^\star\}.
\]
We note that $\EE[\mathbb{I} \{j_\ell\neq i^\star\} \mid \mathcal{F}_{\ell-1}] = 1 - p_{\ell,i_\star}$.
Hoeffding's lemma gives
\[
  \mathbb E\!\left[
    \exp\{\lambda(\mathbb{I} \{j_\ell\neq i^\star\} - (1-p_{\ell,i^\star}))\}
    \mid\mathcal F_{\ell-1}
  \right]\leq e^{\lambda^2/8}.
\]
Applying the tower rule and Markov's inequality, we have
\[
    \PP \left(\exp(\lambda(B_L-P_L)) > \exp\left(\lambda\sqrt{L\log(3/\delta_2)/2}\right)\right)
    \leq \exp\left(\lambda^2 L / 8 - \lambda\sqrt{L\log(3/\delta_2)/2}\right).
\]
Setting $\lambda = 2 \sqrt{2\log(3/\delta_2)/L}$ yields
\begin{equation}
\label{equ:bad-concentation}
    \PP \left(\frac{1}{L}(B_L - P_L) > \sqrt{\frac{\log(3/\delta_2)}{2L}}\right) \le \frac{\delta_2}{3}.
\end{equation}
Because $\beta = 1/(80\bar\alpha) \le 1/80$ and $L = 30\log(3N/\delta_2)/\beta \ge 2400 \log(3/\delta_2)$,
\[
  \sqrt{\frac{\log(3/\delta_2)}{2L}}\leq\frac1{\sqrt{4800}} < \frac{1}{50}.
\]
Combining this with~\eqref{equ:mass-bound}, on the intersection of all three concentration events \eqref{equ:A-bound}, \eqref{equ:mean-concentration}, and \eqref{equ:bad-concentation}
\[
  \frac{B_L}{L} < 0.47 + \frac{1}{50} = 0.49 < \frac12.
\]
Therefore, $B_L < L/2$ with probability at least
$1-\delta_2$.
Hence, $i^\star$ is the unique most frequently deployed controller.

\subsection*{Part 2: System stability}
For this part, we consider the two stages of the algorithm.
In the first $L$ episodes of the algorithm, we are testing different controllers that are potentially destabilizing.
To derive a (conservative) bound for this stage, let us define
\[R_1 := \max\left\{\max_{i, j} \opnorm{A_{[i,j]}}, 2\right\}, R_2 := \max_{i, j} \opnorm{B_{[i,j]}}.\]
And for convenience, denote $T' = \tau L$ and $A_t$ be the closed-loop dynamics at time $t$.
Specifically, we have $A_t = A_{[i_\star, j_{\lceil t / \tau \rceil}]}$.

Then, we directly expand the closed-loop dynamics:
\begin{align*}
\begin{bmatrix} x_{T'} \\ \vdots \\x_2 \\ x_1 \end{bmatrix}
=
\mathcal{T}^{(T')}
\begin{bmatrix} w_{T'} \\ \vdots \\w_2 \\ w_1  \end{bmatrix}
+
\tbo{\mathcal{T}^{(T'-1)}}{\mathbf{0}}
\begin{bmatrix} B_{T'} u_{T'-1} \\ \vdots \\ B_2 u_2 \\ B_1 u_1  \end{bmatrix},
\end{align*}
where
$$\mathcal{T}^{(k)} = 
\begin{bmatrix} 
    I & A_{T'-1} & A_{T'-1}A_{T'-2} & A_{T'-1}A_{T'-2}A_{T'-3} & \dots &A_{T'-1}A_{T'-2} \dots A_{T'-k+1}\\ 
    0 & I & A_{T'-2} & A_{T'-2}A_{T'-3} & \dots & A_{T'-1}A_{T'-2} \dots A_{T'-k+2} \\ 
    0 & 0 & I & A_{T'-3} & \dots & A_{T'-3}A_{T'-4} \dots A_{T'-k+3}\\ 
    \vdots & \vdots & \vdots & \vdots & \ddots & \vdots
\end{bmatrix}$$
is a $k \times k$ block matrix.
In particular, we can write
\[T^{(k)} = \sum_{\ell=0}^{k-1} H^\ell D_\ell,\]
where $H$ is a shift matrix (with 1's on the first super-diagonal and 0's elsewhere), and $D_\ell$ is a block-diagonal matrix corresponding to the diagonals of $\mathcal{T}$.
It follows that 
\[\opnorm{T^{(k)}} \le \sum_{\ell=1}^k R_1^{\ell-1} \le \frac{R_1^{k}}{R_1 -1}.\]
Therefore,
\[\sum_{t=1}^{T'} \norm{x_t}^2 \le 2 \frac{R_1^{2T'}}{(R_1-1)^2} \sum_{t=1}^{T'} \norm{w_t}^2 + 2 \frac{R_1^{2T'-2}R_2^2}{(R_1-1)^2} \sum_{t=1}^{T'} \norm{u_t}^2. \]

Now, when the algorithm succeeds in identifying the matching controller, the closed-loop system follows a stable trajectory after time $T'$.
Specifically, for constants $\kappa_0, \kappa_1$ associated with the matched closed-loop dynamics $A_{[i_\star, i_\star]}$, we have
\[\sum_{t=T'+1}^{T} \norm{x_t}^2 \le \kappa_0 \norm{x_{T'}}^2 + \kappa_1 \sum_{t=T'+1}^{T} (\norm{w_t}^2 + \norm{u_t}^2).\]

In conclusion, we can choose
\[C_0 = 2 (1 + \kappa_0) \max\left\{\frac{R_1^{2T'}}{(R_1-1)^2}, \frac{R_1^{2T'-2}R_2^2}{(R_1-1)^2}\right\}, C_1 = \kappa_1.\]

\end{document}